\documentclass{article}

\usepackage{amssymb}
\usepackage{amsmath}

\usepackage[dvips]{graphicx}

\begin{document}


\begin{center}
\textbf{Multipliers of uniform topological algebras }
\end{center}
 
\noindent \textbf{}

\begin{center}
M. El Azhari 
\end{center}

\noindent \textbf{}

\noindent \textbf{Abstract.} Let $E$ be a complete uniform topological algebra with Arens-Michael normed factors $\left(E_{\alpha}\right)_{\alpha\in\Lambda}.$ Then $M(E) \cong \varprojlim M(E_{\alpha})$ within an algebra isomorphism $\varphi$. If each factor $E_{\alpha}$ is complete, then every multiplier of $E$ is continuous and $\varphi$ is a topological algebra isomorphism where $M(E)$ is endowed with its seminorm topology.
 
\noindent \textbf{}



\noindent \textbf{}
 
\noindent \textbf{1. Preliminaries }                                                                        

\noindent \textbf{}

A topological algebra is an algebra (over the
complex field) which is also a Hausdorff topological vector space such that the
multiplication is separately continuous. For a topological algebra $E,$ we denote by $\Delta\left(E\right)$ the set of all nonzero continuous multiplicative linear functionals on $E.$ An approximate identity in a topological algebra $E$ is a net $ \left(e_{\omega}\right)_{\omega\in\Omega}$ such that for each $x\in E$ we have $xe_{\omega}\rightarrow_{\omega}x$ and $e_{\omega} x\rightarrow_{\omega}x.$ Let $E$ be an algebra, a function $p:E\rightarrow \left[0,\infty\right[$ is called a pseudo-seminorm, if there exists $0\lneq k\leq 1$ such that $p\left(x+y\right)\leq p\left(x\right)+p\left(y\right),\: p\left(\lambda x\right)=\lvert\lambda\rvert^{k}p\left(x\right)$ and $p\left( xy\right)\leq p\left(x\right) p\left(y\right)$ for all $x,y\in E$ and $\lambda\in \mathbb{C}.\; k$ is called the homogenity index of $p.$  If $k=1,\:p$ is called a seminorm. A pseudo-seminorm $p$ is a pseudo-norm, if $p\left(x\right)=0$ implies $x=0.$ A locally \textit{m}-pseudoconvex algebra is a topological algebra $E$ whose topology is determined by a directed family $\left\lbrace p_{\alpha}:\alpha\in\Lambda\right\rbrace$ of pseudo-seminorms. For each $\alpha\in\Lambda,\: ker\left(p_{\alpha}\right) =\left\lbrace x\in E: p_{\alpha}\left(x\right)=0\right\rbrace,$ the quotient algebra $E_{\alpha}=E/ker\left(p_{\alpha}\right)$ is a pseudo-normed algebra in the pseudo-norm $\overline{p}_{\alpha}\left(x_{\alpha}\right)=p_{\alpha}\left(x\right),\: x_{\alpha}=x+ker\left(p_{\alpha}\right).$ Let $f_{\alpha}:E\rightarrow E_{\alpha},\: f_{\alpha}\left(x\right)=x+ker\left(p_{\alpha}\right)=x_{\alpha},$ be the quotient map, $f_{\alpha}$ is a continuous homomorphism from $E$ onto $E_{\alpha}.$ We endow the set $\Lambda$ with the partial order: $\alpha\leq\beta$, if and only, if $p_{\alpha}\left(x\right)\leq p_{\beta}\left(x\right)$ for all $x\in E.$ Take $\alpha\leq\beta$ in $\Lambda,$ since $ker\left(p_{\beta}\right)\subset ker\left(p_{\alpha}\right),$ we define the surjective continuous homomorphism $f_{\alpha\beta}:E_{\beta}\rightarrow E_{\alpha},\: x_{\beta}=x+ker\left(p_{\beta}\right)\rightarrow x_{\alpha}=x+ker\left(p_{\alpha}\right).$ Thus $\left\lbrace\left(E_{\alpha},f_{\alpha\beta}\right),\: \alpha\leq\beta\right\rbrace$ is a projective system of pseudo-normed algebras. We also define the algebra isomorphism (into) $\Phi: E\rightarrow\varprojlim E_{\alpha},\: \Phi(x)=(f_{\alpha}(x))_{\alpha\in\Lambda},$ the canonical projections $\pi_{\alpha}:\prod_{\alpha\in\Lambda} E_{\alpha}\rightarrow E_{\alpha} $ and the restrictions to the projective limit $g_{\alpha}=\pi_{\alpha}/_{\varprojlim E_{\alpha}}:\varprojlim E_{\alpha}\rightarrow E_{\alpha}.$ Since $g_{\alpha}\circ\Phi=f_{\alpha}$ and the quotient map $f_{\alpha}$ is surjective, it follows that the map $ g_{\alpha} $ is surjective, this proves that the projective system  $\left\lbrace\left(E_{\alpha},f_{\alpha\beta}\right),\: \alpha\leq\beta\right\rbrace$ is perfect in the sense of [5, Definition 2.10] (see also [2, Definition 2.7]). Thus, if $E$ is a locally \textit{m}-pseudoconvex algebra (not necessarly complete), then its generalized Arens-Michael projective system  $\left\lbrace\left(E_{\alpha},f_{\alpha\beta}\right),\: \alpha\leq\beta\right\rbrace$ is perfect. If $E$ is complete, then $E\cong\varprojlim E_{\alpha}$ within a topological algebra isomorphism.

\noindent \textbf{} A locally \textit{m}-convex algebra is a topological algebra $E$ whose topology is defined by a directed family $\left\lbrace p_{\alpha}:\alpha\in\Lambda\right\rbrace$ of seminorms. For each $\alpha\in\Lambda,$ put $\Delta_{\alpha}\left(E\right)=\left\lbrace f\in\Delta\left(E\right): \lvert f\left(x\right)\rvert\leq p_{\alpha}\left(x\right),\: x\in E\right\rbrace.$ Let $E$ be an algebra with involution $\ast.$ A seminorm on $E$ is called a $C^{\ast}$-seminorm if $p\left(x^{\ast}x\right)=p\left(x\right)^{2}$ for all $x\in E.$ A complete locally \textit{m}-convex *-algebra $\left(E,\left(p_{\alpha}\right)_{\alpha\in\Lambda}\right),$ for which each $p_{\alpha}$ is a  $C^{\ast}$-seminorm, is called a locally  $C^{\ast}$-algebra. A uniform seminorm on an algebra $E$ is a seminorm $p$ satisfying $p(x^{2})= p(x)^{2}$ for all $x\in E.$ A uniform topological algebra is a topological algebra whose topology is determined by a directed family of
uniform seminorms. In that case, such a topological algebra is also named a \textit{uniform locally convex algebra}. A uniform normed algebra is a normed algebra $(E,\Vert .\Vert)$ such
that  $\Vert x^{2}\Vert = \Vert x\Vert^{2}$ for all $x\in E.$  

\noindent \textbf{} An algebra $E$ is called proper if for any $x\in E,\: xE=Ex=\left\lbrace 0\right\rbrace$ implies $x=0.$ If $E$ has identity, then $E$ is proper. Moreover, a topological algebra with approximate identity is proper. Also, a (Hausdorff) uniform topological algebra is proper. Let $E$ be an algebra, a map $T:E\rightarrow E$ is called a multiplier if $T\left(x\right)y=xT\left(y\right)$ for all $x,y \in E.$ We denote by $M(E)$ the set of all multipliers of $E.$ It is known that if $E$ is a proper algebra, then any multiplier $T$ of $E$ is linear with the property  $T\left(xy\right)=T\left(x\right)y=xT\left(y\right)$ for all $x,y \in E,$ and $M(E)$ is a commutative algebra with the identity map $ I $ of $E$ as its identity. Let $\left(E,\lVert .\rVert\right)$ be a uniform normed algebra, and let $M_{c}(E)$ be the algebra of all continuous multipliers of $E$ with the operator norm $\lVert .\rVert_{op}.$ It is known that $\lVert .\rVert_{op}$ has the square property and the map $l:\left(E,\lVert .\rVert\right)\rightarrow\left(M_{c}(E),\lVert .\rVert_{op}\right),\: l\left(x\right)\left(y\right)=xy,$ is an isometric isomorphism (into). For information on the multiplier algebra in non-normed topological algebras, see also [3] and [4].

\noindent \textbf{} In the sequel, we will need the following elementary result called the universal property of the quotient: Let $ X,Y,Z $ be vector spaces, $ f: X\rightarrow Y $ and $ g: X\rightarrow Z $ be linear maps. If the map $g$ is surjective and $ker(g)\subset ker(f),$ then there exists a unique linear map $ h: Z\rightarrow Y $ such that $ f=h\circ g. $

\noindent \textbf{}

\noindent \textbf{2. Results}

\noindent \textbf{}

\noindent \textbf{Proposition 2.1.} \textit{Let $\left(E,\left(p_{\alpha}\right)_{\alpha\in\Lambda}\right)$ be a locally \textit{m}-pseudoconvex algebra with proper pseudo-normed factors $ \left( E_{\alpha}\right)_{\alpha\in\Lambda}.$ The following assertions are equivalent:\\
(i) $T\left(ker\left(f_{\alpha}\right) \right)\subset ker\left(f_{\alpha}\right)$ for all $T\in M(E)$ and $\alpha\in\Lambda;$\\
(ii) for each $T\in M(E),$ there exists a unique $\left(T_{\alpha}\right)_{\alpha\in\Lambda}\in\prod_{\alpha\in\Lambda} M(E_{\alpha})$ such that $f_{\alpha}\circ T=T_{\alpha}\circ f_{\alpha}$ and $T_{\alpha}\circ f_{\alpha\beta}=f_{\alpha\beta}\circ T_{\beta}$ for all $\alpha\leq\beta$ in $\Lambda.$ Furthermore, $T$ is continuous if and only if $T_{\alpha}$ is continuous for all $\alpha\in\Lambda.$}

\noindent \textbf{}
 
\noindent \textbf {Proof.} Since the pseudo-normed factors $(E_{\alpha})_{\alpha\in\Lambda}$ are proper, it follows that the algebra $E$ is proper and so every multiplier of $E$ (or $E_{\alpha}$ ) is linear.\\
$\left(ii\right)\Rightarrow\left(i\right):$ If $T\in M(E)$ and $x\in ker\left(f_{\alpha}\right),$ then $f_{\alpha}\left(T\left(x\right)\right)=T_{\alpha}\left(f_{\alpha}\left(x\right)\right)=0$ and so $T\left(x\right)\in ker\left(f_{\alpha}\right) .$\\
$\left(i\right)\Rightarrow\left(ii\right):$ Take $T\in M(E)$ and $\alpha\in\Lambda.$ Since $T\left(ker\left(f_{\alpha}\right) \right)\subset ker\left(f_{\alpha}\right)$ and by using the universal property of the quotient (see Preliminaries), there exists a unique linear map $T_{\alpha}:E_{\alpha}\rightarrow E_{\alpha}$ such that $f_{\alpha}\circ T=T_{\alpha}\circ f_{\alpha}.$ Let $\alpha\in\Lambda$ and $x,y\in E,
\;T_{\alpha}\left(f_{\alpha}\left(x\right)f_{\alpha}\left(y\right)\right)=T_{\alpha}\left(f_{\alpha}\left(xy\right)\right)=f_{\alpha}\left(T\left(xy\right)\right)=f_{\alpha}\left(xT\left(y\right)\right)=f_{\alpha}\left(x\right)f_{\alpha}\left(T\left(y\right)\right)=f_{\alpha}\left(x\right)T_{\alpha}\left(f_{\alpha}\left(y\right)\right)$ and similary on the other side, so $T_{\alpha}$ is a multiplier of $E_{\alpha}.$ Let $\alpha\leq\beta$ in $\Lambda,$ we have $T_{\alpha}\circ f_{\alpha}=f_{\alpha}\circ T,$ then $T_{\alpha}\circ f_{\alpha\beta}\circ f_{\beta}=f_{\alpha\beta}\circ f_{\beta}\circ T=f_{\alpha\beta}\circ T_{\beta}\circ f_{\beta},$ hence  $T_{\alpha}\circ f_{\alpha\beta}=f_{\alpha\beta}\circ T_{\beta}$ since the quotient map $f_{\beta}$ is surjective. Suppose that $T$ is continuous. Let $O_{\alpha}$ be an open set in $E_{\alpha},$ we have $f_{\alpha}^{-1}\left(T_{\alpha}^{-1}\left(O_{\alpha}\right)\right)=\left(T_{\alpha}\circ f_{\alpha}\right)^{-1}\left(O_{\alpha}\right)=\left(f_{\alpha}\circ T\right)^{-1}\left(O_{\alpha}\right)$ which is open in $E$ since $f_{\alpha}\circ T$ is continuous, then $T_{\alpha}^{-1}\left(O_{\alpha}\right)$ is open in $E_{\alpha}.$ Conversely, suppose that $T_{\alpha}$ is continuous for all $\alpha\in\Lambda.$ Since $E$ is topologically isomorphic to a subalgebra of $\varprojlim E_{\alpha},\; T$ is continuous if and only if $f_{\alpha}\circ T$ is continuous for all $\alpha\in\Lambda.$ Since $f_{\alpha}\circ T=T_{\alpha}\circ f_{\alpha}$ and $T_{\alpha}$ is continuous for all $\alpha\in\Lambda,$ we deduce that $T$ is continuous.
 
\noindent \textbf{}

\noindent \textbf{Proposition 2.2.} \textit{Let $\left(E,\left(p_{\alpha}\right)_{\alpha\in\Lambda}\right)$ be a locally \textit{m}-pseudoconvex algebra with proper pseudo-normed factors $ \left( E_{\alpha}\right)_{\alpha\in\Lambda}.$ The following assertions are equivalent:\\ 
(j) $U\left(ker\left(f_{\alpha\beta}\right)\right)\subset ker\left(f_{\alpha\beta}\right)$ for all $U\in M(E_{\beta})$ and $\alpha\leq\beta$ in $\Lambda;$\\
(jj) there exists a unique projective system $\left\lbrace\left(M(E_{\alpha}),h_{\alpha\beta}\right),\:\alpha\leq\beta\right\rbrace$ such that $h_{\alpha\beta}\left(U\right)\circ f_{\alpha\beta}=f_{\alpha\beta}\circ U$ for all $U\in M(E_{\beta})$ and $\alpha\leq\beta$ in $\Lambda.$ Furthermore, if $E_{\alpha}$ is complete for all $\alpha\in\Lambda,$ then $h_{\alpha\beta}$ is continuous for all $\alpha\leq\beta$ in $\Lambda.$}

\noindent \textbf{}
 
\noindent \textbf{Proof.} $\left(jj\right)\Rightarrow\left(j\right):$ Let $U\in M(E_{\beta})$ and $x_{\beta}\in ker\left(f_{\alpha\beta}\right),$ then $f_{\alpha\beta}\left(U\left(x_{\beta}\right)\right)=h_{\alpha\beta}\left(U\right)\left(f_{\alpha\beta}\left(x_{\beta}\right)\right)=0$ and so $U\left(x_{\beta}\right)\in ker\left(f_{\alpha\beta}\right).$\\
$\left(j\right)\Rightarrow\left(jj\right):$ Let $\alpha\leq\beta$ in $\Lambda$ and $U\in M(E_{\beta}).$ Since $U\left(ker\left(f_{\alpha\beta}\right)\right)\subset ker\left(f_{\alpha\beta}\right)$ and by using the universal property of the quotient (see Preliminaries), there exists a unique linear map $V:E_{\alpha}\rightarrow E_{\alpha}$ such that $V\circ f_{\alpha\beta}=f_{\alpha\beta}\circ U.$ Let $x_{\alpha}=x+ ker(p_{\alpha}), y_{\alpha}=y+ ker(p_{\alpha})\in E_{\alpha}$ where $x,y\in E.$ Put $x_{\beta}=x+ ker(p_{\beta})$ and $ y_{\beta}=y+ ker(p_{\beta}),$ clearly $ x_{\beta}, y_{\beta}\in E_{\beta}.$ By definition of the map $f_{\alpha\beta},$ we get $ f_{\alpha\beta}(x_{\beta})= x_{\alpha}$ and $f_{\alpha\beta}(y_{\beta})= y_{\alpha}.$ We have  $V(x_{\alpha}y_{\alpha})= V\left(f_{\alpha\beta}\left( x_{\beta}\right) f_{\alpha\beta}\left( y_{\beta}\right) \right) =V\left( f_{\alpha\beta}\left( x_{\beta}y_{\beta}\right) \right) =f_{\alpha\beta}\left( U\left( x_{\beta}y_{\beta}\right) \right) =f_{\alpha\beta}\left( x_{\beta}U\left( y_{\beta}\right) \right) =f_{\alpha\beta}\left( x_{\beta}\right) f_{\alpha\beta}\left( U\left( y_{\beta}\right) \right) =f_{\alpha\beta}\left( x_{\beta}\right) V\left( f_{\alpha\beta}\left( y_{\beta}\right) \right)= x_{\alpha}V(y_{\alpha})$ and similary on the other side, so $V$ is a multiplier of $E_{\alpha}.$ This shows the existence of the map $ h_{\alpha\beta}:M( E_{\beta}) \rightarrow M( E_{\alpha})$ such that $h_{\alpha\beta}\left(U\right)\circ f_{\alpha\beta}=f_{\alpha\beta}\circ U$ for all $U\in M(E_{\beta})$ and $\alpha\leq\beta$ in $\Lambda.$ Let $ \alpha\leq\beta $ in $ \Lambda,\:U_{1},U_{2}\in M( E_{\beta}) $ and $ \lambda\in\mathbb{C},\;h_{\alpha\beta}\left( U_{1}+\lambda U_{2}\right) \circ f_{\alpha\beta}=f_{\alpha\beta}\circ\left( U_{1}+\lambda U_{2}\right) =\left( f_{\alpha\beta}\circ U_{1}\right) +\lambda\left( f_{\alpha\beta}\circ U_{2}\right) =h_{\alpha\beta}\left( U_{1}\right) \circ f_{\alpha\beta}+\lambda h_{\alpha\beta}\left( U_{2}\right) \circ f_{\alpha\beta} =\left( h_{\alpha\beta}\left( U_{1}\right) +\lambda h_{\alpha\beta}\left( U_{2}\right) \right) \circ f_{\alpha\beta},$ hence\\
$ h_{\alpha\beta}\left( U_{1}+\lambda U_{2}\right) =h_{\alpha\beta}\left( U_{1}\right) +\lambda h_{\alpha\beta}\left( U_{2}\right) $ since $ f_{\alpha\beta} $ is surjective. Also,\\
$ h_{\alpha\beta}\left( U_{1}\circ U_{2}\right) \circ f_{\alpha\beta}=f_{\alpha\beta}\circ U_{1}\circ U_{2} =h_{\alpha\beta}\left( U_{1}\right) \circ f_{\alpha\beta}\circ U_{2}=h_{\alpha\beta}\left( U_{1}\right) \circ h_{\alpha\beta}\left( U_{2}\right) \circ f_{\alpha\beta},$ then $ h_{\alpha\beta}\left( U_{1}\circ U_{2}\right) =h_{\alpha\beta}\left( U_{1}\right) \circ h_{\alpha\beta}\left( U_{2}\right)$ since $ f_{\alpha\beta} $ is surjective. Let $ \alpha\leq\beta\leq\gamma $ in $ \Lambda $ and $ W\in M( E_{\gamma}),\;\left( h_{\alpha\beta}\circ h_{\beta\gamma}\right) \left( W\right) \circ f_{\alpha\gamma}=h_{\alpha\beta}\left( h_{\beta\gamma}\left( W\right) \right) \circ f_{\alpha\beta} \circ f_{\beta\gamma}=f_{\alpha\beta}\circ h_{\beta\gamma}\left( W\right) \circ f_{\beta\gamma}=f_{\alpha\beta}\circ f_{\beta\gamma}\circ W=f_{\alpha\gamma}\circ W=h_{\alpha\gamma}\left( W\right) \circ f_{\alpha\gamma},$ consequently $ \left( h_{\alpha\beta}\circ h_{\beta\gamma}\right) \left( W\right) =h_{\alpha\gamma}\left( W\right) $ since $ f_{\alpha\gamma} $ is surjective. Thus $ h_{\alpha\beta}\circ h_{\beta\gamma}=h_{\alpha\gamma}. $ Let $ \alpha\in\Lambda, $ if $ E_{\alpha} $ is complete, then every multiplier of $ E_{\alpha} $ is continuous. Now by assuming that $ E_{\alpha} $ is complete for all $ \alpha\in\Lambda, $ we will show that $ h_{\alpha\beta} $ is continuous for all $ \alpha\leq\beta $ in $ \Lambda$ (see also, the proof of Theorem 2.12 in [5]). For $ \alpha\in\Lambda$ and $r\gneq 0,  $ let $ B_{\alpha}\left( 0,r\right)=\left\lbrace x_{\alpha}\in E_{\alpha}:\overline{p}_{\alpha}\left( x_{\alpha}\right) \leq r\right\rbrace.  $  We denote by $ \lVert .\rVert_{\alpha} $ the operator pseudo-norm on $ M\left( E_{\alpha} \right). $ Let $ \alpha\leq\beta $ in $ \Lambda,\:f_{\alpha\beta} $ is open by the open mapping theorem, so there is $ \lambda\gneq 0 $ such that $ \lambda B_{\alpha}\left( 0,1\right) \subset f_{\alpha\beta}\left( B_{\beta}\left( 0,1\right) \right)  $ i.e. $ B_{\alpha}\left( 0,1\right) \subset f_{\alpha\beta}\left( B_{\beta} \left( 0,r\right) \right) $ where $ r=\lambda^{-k_{\beta}} $ and $ k_{\beta} $ is the homogenity index of $ \overline{p}_{\beta}. $ Let $ U\in M\left( E_{\beta} \right),\\
\lVert h_{\alpha\beta}\left( U\right) \rVert_{\alpha}=\sup \left\lbrace \overline{p}_{\alpha}\left( h_{\alpha\beta}\left( U\right) \left( f_{\alpha}\left( x\right) \right)\right)  :  f_{\alpha}\left( x\right) \in B_{\alpha}\left( 0,1\right) \right\rbrace \\
\leq \sup \left\lbrace \overline{p}_{\alpha}\left(h_{\alpha\beta}\left( U\right) \left( f_{\alpha\beta}\left( f_{\beta}\left( x\right) \right) \right) \right) :f_{\beta}\left( x\right) \in B_{\beta}\left( 0,r\right) \right\rbrace \\
=\sup \left\lbrace \overline{p}_{\alpha}\left( f_{\alpha\beta}\left( U\left( f_{\beta}\left( x\right) \right) \right) \right) :f_{\beta}\left( x\right) \in B_{\beta}\left( 0,r\right) \right\rbrace \\
\leq\sup \left\lbrace \overline{p}_{\beta}\left( U\left( f_{\beta}\left( x\right) \right) \right) :f_{\beta}\left( x\right) \in B_{\beta}\left( 0,r\right) \right\rbrace \\
\leq\sup \left\lbrace \lVert U\rVert_{\beta}\:\overline{p}_{\beta}\left( f_{\beta}\left( x\right) \right) :f_{\beta}\left( x\right) \in B_{\beta}\left( 0,r\right) \right\rbrace \\
= r\lVert U\rVert_{\beta}. $  \\
Therefore $ h_{\alpha\beta}$ is continuous.

\noindent \textbf{}

\noindent \textbf{Theorem 2.3.} \textit{Let $\left(E,\left(p_{\alpha}\right)_{\alpha\in\Lambda}\right)$ be a complete locally \textit{m}-pseudoconvex algebra with proper pseudo-normed factors $ \left( E_{\alpha}\right)_{\alpha\in\Lambda}.$ Assume that $ E $ satisfies conditions (i) and (j). Then $M(E) \cong \varprojlim M(E_{\alpha})$ within an algebra isomorphism $\varphi$. Furthermore, if each factor $ E_{\alpha} $ is complete, then every multiplier of $ E $ is continuous and $ \varphi $ is a topological algebra isomorphism where $M(E)$ is endowed with its pseudo-seminorm topology.}
 
\noindent \textbf{}
 
\noindent \textbf{Proof.} Take $T\in M(E).$ By Propositions 2.1 and 2.2, $\left(T_{\alpha}\right)_{\alpha\in\Lambda}\in\prod_{\alpha\in\Lambda} M(E_{\alpha}),\\
T_{\alpha}\circ f_{\alpha\beta}=f_{\alpha\beta}\circ T_{\beta}$ and $h_{\alpha\beta}\left(T_{\beta}\right)\circ f_{\alpha\beta}=f_{\alpha\beta}\circ T_{\beta}$ for all $\alpha\leq\beta$ in $\Lambda.$ Hence $ h_{\alpha\beta}(T_{\beta})\circ f_{\alpha\beta}= T_{\alpha}\circ f_{\alpha\beta}$ and consequently $h_{\alpha\beta}(T_{\beta})= T_{\alpha}$ since the map $f_{\alpha\beta}$ is surjective. This shows that $(T_{\alpha})_{\alpha\in\Lambda}\in\varprojlim M(E_{\alpha}).$ Thus the map $ \varphi:M(E)\rightarrow \varprojlim M(E_{\alpha}),\: 
T\rightarrow \left( T_{\alpha}\right) _{\alpha\in\Lambda},$ is well defined. We will show that $ \varphi $ is an algebra isomorphism. Let $ T,S\in M\left( E\right)  $ and $ \lambda\in\mathbb{C},\;T_{\alpha}\circ f_{\alpha}=f_{\alpha}\circ T $ and $ S_{\alpha}\circ f_{\alpha}=f_{\alpha}\circ S, $ then $ \left( T_{\alpha}+\lambda S_{\alpha}\right) \circ f_{\alpha}=f_{\alpha} \circ\left( T+\lambda S\right),$ so $ \left( T+\lambda S\right)_{\alpha} =T_{\alpha}+\lambda S_{\alpha}$ by Proposition 2.1. Also, $ T_{\alpha}\circ S_{\alpha}\circ f_{\alpha}=T_{\alpha}\circ f_{\alpha} \circ S=f_{\alpha}\circ T\circ S,$ hence $ \left( T\circ S\right)_{\alpha}=T_{\alpha} \circ S_{\alpha} $ by Proposition 2.1. Let $ T\in M\left( E\right), $ if $ T_{\alpha}=0 $ for all $ \alpha\in\Lambda, $ then $ f_{\alpha}\circ T=T_{\alpha}\circ f_{\alpha}=0 $ for all $ \alpha\in\Lambda $ and consequently $ T=0. $ Let $ \left( U_{\alpha}\right)_{\alpha\in\Lambda} \in \varprojlim M(E_{\alpha})$ and define the map $ T=\Phi^{-1}\circ\varprojlim U_{\alpha}\circ \Phi:E\rightarrow E $ where $\varprojlim U_{\alpha}$ is the multiplier of $ \varprojlim E_{\alpha} $ defined by $  \left( \varprojlim U_{\alpha}\right) \left( x_{\alpha}\right)_{\alpha}=\left( U_{\alpha}\left( x_{\alpha} \right) \right)_{\alpha}$ and $ \Phi:E\rightarrow \varprojlim E_{\alpha} $ is the topological algebra isomorphism given by $ \Phi\left( x\right) =\left( f_{\alpha}\left( x\right) \right)_{\alpha}. $ Clearly $ T $ is a multiplier of $ E, $ also $ f_{\alpha}\circ T=f_{\alpha}\circ\Phi^{-1}\circ\varprojlim U_{\alpha}\circ \Phi=U_{\alpha}\circ f_{\alpha}$ for all $ \alpha\in\Lambda, $ so $ \varphi\left( T\right) =\left( U_{\alpha}\right)_{\alpha}.$ If $ E_{\alpha} $ is complete for all $ \alpha\in\Lambda, $ then every multiplier of $ E_{\alpha} $ is continuous, hence every multiplier of $ E $ is continuous by Proposition 2.1. The pseudo-seminorm topology on $M(E)$ is the topology defined by the family of pseudo-seminorms $ q_{\alpha}\left( T\right) =\lVert T_{\alpha}\rVert_{\alpha},\:\alpha\in\Lambda, $ so $ \varphi $ is a topological algebra isomorphism.

\noindent \textbf{}
 
\noindent \textbf{Proposition 2.4.} \textit{Let $\left(E,\left(p_{\alpha}\right)_{\alpha\in\Lambda}\right)$ be a locally \textit{m}-pseudoconvex algebra with approximate identity $ \left( e_{\omega}\right)_{\omega\in\Omega}.$ Then $ E $ satisfies conditions (i) and (j).}
 
\noindent \textbf{}

\noindent \textbf{Proof.} Let $ T\in M\left( E\right),\: x\in ker\left( f_{\alpha}\right) $ and $ \omega\in\Omega, $  \\
$ f_{\alpha}\left( T\left( x\right) \right) =f_{\alpha}\left( T\left( x-xe_{\omega}+xe_{\omega}\right) \right) \\
= f_{\alpha} \left( T\left( x\right) -T\left( xe_{\omega}\right) \right) +f_{\alpha}\left( T\left( xe_{\omega}\right) \right) =f_{\alpha}\left( T\left( x\right) -T\left( x\right) e_{\omega}\right)+ f_{\alpha}\left( xT\left( e_{\omega}\right) \right) \\
= f_{\alpha}\left( T\left( x\right) -T\left( x\right) e_{\omega}\right)+f_{\alpha}\left( x\right) f_{\alpha}\left( T\left( e_{\omega}\right) \right) =f_{\alpha}\left( T\left( x\right) -T\left( x\right) e_{\omega}\right).$\\
Since $ T\left( x\right) e_{\omega}\rightarrow_{\omega} T\left( x\right)  $ and $ f_{\alpha} $ is continuous, we deduce that $ f_{\alpha}\left( T\left( x\right) \right) =0. $ Now we will show that $ U\left( ker\left( f_{\alpha\beta}\right) \right) \subset ker\left( f_{\alpha\beta}\right)   $ for all $ U\in M\left( E_{\beta}\right)  $
and $ \alpha\leq\beta $ in $ \Lambda. $ Since $ \left( e_{\omega}\right)_{\omega\in\Omega}$ is an approximate identity in $ E $ and $ f_{\beta}:E\rightarrow E_{\beta} $ is a surjective continuous homomorphism, it follows that $ \left( f_{\beta}\left( e_{\omega}\right)\right) _{\omega\in\Omega}$ is an approximate identity in $ E_{\beta}$ (see [8, Theorem 4.1]). Let $ U\in M\left( E_{\beta} \right),\: x_{\beta}\in ker\left( f_{\alpha\beta}\right)$ and  $ \omega\in\Omega, $\\
$ f_{\alpha\beta}\left( U\left( x_{\beta}\right) \right) =f_{\alpha\beta}\left( U\left( x_{\beta}-x_{\beta}f_{\beta}\left( e_{\omega}\right) + x_{\beta}f_{\beta}\left( e_{\omega}\right)\right) \right) \\
=f_{\alpha\beta}\left( U\left( x_{\beta}\right) -U\left( x_{\beta}f_{\beta}\left( e_{\omega}\right) \right) \right) + f_{\alpha\beta}\left( U\left( x_{\beta}f_{\beta}\left( e_{\omega}\right) \right) \right) \\
=f_{\alpha\beta}\left( U\left( x_{\beta}\right) -U\left( x_{\beta}\right) f_{\beta}\left( e_{\omega}\right) \right) + f_{\alpha\beta}\left( x_{\beta}U\left( f_{\beta}\left( e_{\omega}\right) \right) \right) \\
= f_{\alpha\beta}\left( U\left( x_{\beta}\right) -U\left( x_{\beta}\right) f_{\beta}\left( e_{\omega}\right) \right) + f_{\alpha\beta}\left( x_{\beta}\right) f_{\alpha\beta}\left( U\left( f_{\beta}\left( e_{\omega}\right) \right) \right)\\
= f_{\alpha\beta}\left( U\left( x_{\beta}\right) -U\left( x_{\beta}\right) f_{\beta}\left( e_{\omega}\right) \right).$\\
Since $ U\left( x_{\beta}\right) f_{\beta}\left( e_{\omega}\right) \rightarrow_{\omega} U\left( x_{\beta}\right) $ and $ f_{\alpha\beta} $ is continuous, we deduce that\\
$ f_{\alpha\beta}\left( U\left( x_{\beta}\right) \right)=0.  $

\noindent \textbf{}

\noindent \textbf{Corollary 2.5.} [5, Theorems 2.6 and 2.12] \textit{Let $\left(E,\left(p_{\alpha}\right)_{\alpha\in\Lambda}\right)$ be a complete locally m-pseudoconvex algebra with approximate identity. Suppose that each factor $ E_{\alpha}=E/ker\left( p_{\alpha}\right) $ in the generalized Arens-Michael decomposition of $ E $ is complete. Then every multiplier of $ E $ is continuous and $M(E)\cong \varprojlim M(E_{\alpha})$ within a topological algebra isomorphism where $M(E)$ is endowed with its pseudo-seminorm topology.}
 
\noindent \textbf{}

\noindent \textbf{Proof.} It follows from Theorem 2.3 and Proposition 2.4.
 
\noindent \textbf{}

\noindent \textbf{Corollary 2.6.} \textit{Let $\left(E,\left(p_{\alpha}\right)_{\alpha\in\Lambda}\right)$ be a locally $C^{\ast}$-algebra. Then every multiplier of $ E $ is continuous and $M(E)\cong\varprojlim M(E_{\alpha})$ within a topological algebra isomorphism where $M(E)$ is endowed with its seminorm topology.}

\noindent \textbf{}

\noindent \textbf{Proof.} By [7, Theorem 2.6] and [10, Corollary 1.12], $ E $ has an approximate identity and each factor $ E_{\alpha} $ is complete.
 
\noindent \textbf{}

Now we will describe multiplier algebras of complete uniform topological algebras.
 
\noindent \textbf{} 

\noindent \textbf{Proposition 2.7.} \textit{Let $\left(E,\left(p_{\alpha}\right)_{\alpha\in\Lambda}\right)$ be a uniform topological algebra. Then\\
$ ker\left( f_{\alpha}\right) =\cap\left\lbrace ker\left( \chi\right) :\chi\in\Delta_{\alpha}\left( E\right) \right\rbrace  $ for all $ \alpha\in\Lambda $ and\\
$ ker\left( f_{\alpha\beta}\right) =\cap\left\lbrace ker\left( \mu\circ f_{\alpha\beta}\right) :\mu\in \Delta\left( E_{\alpha} \right) \right\rbrace $ for all $ \alpha\leq\beta $ in $ \Lambda. $}
 
\noindent \textbf{}

\noindent \textbf{Proof.} Show first that $\Delta_{\alpha}(E)$ and $ \Delta(E_{\alpha})$ are non empty sets. Let $ F_{\alpha} $ be the completion of $ (E_{\alpha},\overline{p}_{\alpha}),\: F_{\alpha}$ is a uniform Banach algebra. By [8, Lemma 5.1], $F_{\alpha}$ is commutative and semisimple. Then $ \Delta(F_{\alpha})$ is a non empty set since $F_{\alpha}$ is not a radical algebra, hence $\Delta_{\alpha}(E)$ and $ \Delta(E_{\alpha})$ are non empty sets (see [9, Proposition 7.5]).\\
By [1, Theorem 6], $ p_{\alpha}\left( x\right) =\sup \left\lbrace \lvert\chi\left( x\right) \rvert : \chi\in\Delta_{\alpha}\left( E\right) \right\rbrace  $ for all $ x\in E $ and $ \alpha\in\Lambda, $ then $ ker\left( f_{\alpha}\right)= ker\left( p_{\alpha}\right) =\cap\left\lbrace ker\left( \chi\right) :\chi\in\Delta_{\alpha}\left( E\right) \right\rbrace  $ for all $ \alpha\in\Lambda. $ Let $ \alpha\leq\beta $ in $ \Lambda $ and $ x_{\beta}\in E_{\beta}$  ,\\
$ x_{\beta}\in ker\left( f_{\alpha\beta}\right) \Leftrightarrow f_{\alpha\beta}\left( x_{\beta}\right) =0$\\
$ \Leftrightarrow\mu\left( f_{\alpha\beta}\left( x_{\beta}\right) \right) =0 $ for all $ \mu\in\Delta\left( E_{\alpha}\right)$\\
$ \Leftrightarrow x_{\beta}\in \cap\left\lbrace ker\left( \mu\circ f_{\alpha\beta}\right) : \mu\in\Delta\left( E_{\alpha}\right) \right\rbrace.$

\noindent \textbf{} 

\noindent \textbf{Proposition 2.8.} \textit{Let $\left(E,\left(p_{\alpha}\right)_{\alpha\in\Lambda}\right)$ be a uniform topological algebra. Then $ E $ satisfies conditions (i) and (j).}
 
\noindent \textbf{} 

\noindent \textbf{Proof.} By Proposition 2.7, $ ker\left( f_{\alpha}\right) =\cap\left\lbrace ker\left( \chi\right) :\chi\in\Delta_{\alpha}\left( E\right) \right\rbrace  $ for all $ \alpha\in\Lambda. $ If $ T $ is a multiplier of $ E, $ then $ T\left( ker\left( \chi\right) \right) \subset ker\left( \chi\right)  $ for all $ \chi\in\Delta_{\alpha}\left( E\right)  $ by [6, Theorem 2.9] and [8, Lemma 5.1], so\\
$ T\left( ker\left( f_{\alpha}\right) \right) =T\left( \cap\left\lbrace ker\left( \chi\right) :\chi\in\Delta_{\alpha}\left( E\right) \right\rbrace \right)  $\\
$ \subset\cap\left\lbrace T\left( ker\left( \chi\right) \right) :\chi\in\Delta_{\alpha}\left( E\right) \right\rbrace  $\\
$ \subset\cap\left\lbrace ker\left( \chi\right) :\chi\in\Delta_{\alpha}\left( E\right) \right\rbrace = ker\left( f_{\alpha}\right). $\\
By Proposition 2.7, $ ker\left( f_{\alpha\beta}\right) =\cap\left\lbrace ker\left( \mu\circ f_{\alpha\beta}\right) :\mu\in \Delta\left( E_{\alpha} \right) \right\rbrace $ for all $ \alpha\leq\beta $ in $ \Lambda. $ If $ U $ is a multiplier of $ E_{\beta}, $ then $ U\left( ker\left( \delta\right) \right) \subset ker\left( \delta\right)  $ for all $ \delta\in\Delta\left( E_{\beta}\right)  $ by [6, Theorem 2.9] and [8, Lemma 5.1], so $ U\left( ker\left( \mu\circ f_{\alpha\beta}\right) \right) \subset ker\left( \mu\circ f_{\alpha\beta}\right) $ for all $ \mu\in\Delta\left( E_{\alpha}\right), $ and consequently\\
$ U\left( ker\left( f_{\alpha\beta}\right) \right) = U\left( \cap\left\lbrace ker\left( \mu\circ f_{\alpha\beta}\right) : \mu\in\Delta\left( E_{\alpha}\right) \right\rbrace \right)  $\\
$ \subset\cap\left\lbrace U\left( ker\left( \mu\circ f_{\alpha\beta}\right) \right) :\mu\in\Delta\left( E_{\alpha} \right) \right\rbrace  $\\
$ \subset\cap\left\lbrace ker\left( \mu\circ f_{\alpha\beta} \right) : \mu\in\Delta\left( E_{\alpha}\right) \right\rbrace = ker\left( f_{\alpha\beta}\right).$
 
\noindent \textbf{}  

\noindent \textbf{Theorem 2.9.} \textit{Let  $\left(E,\left(p_{\alpha}\right)_{\alpha\in\Lambda}\right)$ be a complete uniform topological algebra. Then $M(E) \cong \varprojlim M(E_{\alpha})$ within an algebra isomorphism $\varphi$. Furthermore, if each factor $E_{\alpha}$ is complete, then every multiplier of $E$ is continuous and $\varphi$ is a topological algebra isomorphism where $M(E)$ is endowed with its seminorm topology.}

\noindent \textbf{} 

\noindent \textbf{Proof.} It follows from Theorem 2.3 and Proposition 2.8.

\noindent \textbf{} 

\noindent \textbf{Remark.} Let  $\left(E,\left(p_{\alpha}\right)_{\alpha\in\Lambda}\right)$ be a complete uniform topological algebra which is also a symmetric *-algebra. Then $ \chi(x^{\ast})=\overline{\chi(x)}$ for all $x\in E$ and $ \chi\in\Delta(E)$ (see [9, Lemma 6.4]). Take $ x\in E $ and $ \alpha\in\Lambda.$ By [1, Theorem 6],\\
$ p_{\alpha}(x^{\ast} x)=\sup\left\lbrace\lvert\chi(x^{\ast} x)\rvert:\;\chi\in\Delta_{\alpha}(E)\right\rbrace 
=\sup\left\lbrace \lvert\chi(x)\rvert^{2}:\;\chi\in\Delta_{\alpha}(E)\right\rbrace \\
=(\sup\left\lbrace \lvert\chi(x)\rvert :\;\chi\in\Delta_{\alpha}(E)\right\rbrace )^{2}
= p_{\alpha}(x)^{2}.$\\
Therefore $\left(E,\left(p_{\alpha}\right)_{\alpha\in\Lambda}\right)$ is a locally  $C^{\ast}$-algebra, and so each factor $ E_{\alpha} $ is complete.

\noindent \textbf{} 

As an application of previous results, we deduce the Arhippainen unitization theorem [1, Theorem 4] on uniform topological algebras.

\noindent \textbf{} 

\noindent \textbf{Proposition 2.10.} \textit{Let  $\left(E,\left(p_{\alpha}\right)_{\alpha\in\Lambda}\right)$ be a uniform topological algebra, and let $ M_{c}( E)  $ be the algebra of all continuous multipliers of $ E. $ Then there is a family of seminorms $ \left( q_{\alpha}\right)_{\alpha\in\Lambda} $ on $ M_{c}( E)  $ such that\\
1. $ \left( M_{c}(E) ,\left( q_{\alpha}\right)_{\alpha\in\Lambda}\right)  $ is a uniform topological algebra;\\
2. the map $ L:E\rightarrow M_{c}(E),\:L\left( x\right) \left( y\right) =xy, $ is an algebra isomorphism (into) and $ q_{\alpha}\left( L\left( x\right) \right) =p_{\alpha}\left( x\right)  $ for all $ x\in E $ and $ \alpha\in\Lambda. $}

\noindent \textbf{} 

\noindent \textbf{Proof.} 1. By Propositions 2.1, 2.2 and 2.8, we define the map  $ \psi:M_{c}(E) \rightarrow \varprojlim M_{c}(E_{\alpha}),\: T\rightarrow \left( T_{\alpha}\right) _{\alpha\in\Lambda}. $ As in the proof of Theorem 2.3, $ \psi $ is an injective homomorphism. We endow $ M_{c}(E)  $ with the topology defined by the family of seminorms $ q_{\alpha}\left( T\right) =\lVert T_{\alpha}\rVert_{\alpha},\:\alpha\in\Lambda,$ where $ \lVert .\rVert_{\alpha} $ is the operator norm on $ M_{c}( E_{\alpha}).$ Let $ T\in M_{c}(E) ,\:q_{\alpha}\left( T^{2}\right) =\lVert \left( T^{2}\right)_{\alpha}\rVert_{\alpha} =\lVert\left( T_{\alpha}\right)^{2}\rVert_{\alpha} =\lVert T_{\alpha}\rVert_{\alpha}^{2}= q_{\alpha}\left( T\right)^{2}$ since $ \lVert .\rVert_{\alpha} $ has the square property. Let $ T\in M_{c}(E)  $ with $ q_{\alpha}\left( T\right) =0 $ for all $ \alpha\in\Lambda, $ then $ T_{\alpha}=0 $ for all $ \alpha\in\Lambda, $ so $ T=0$  since $ \psi $ is injective.\\
2. Since $ E $ is proper, $ L $ is an algebra isomorphism (into). Let $ x\in E $ and $ \alpha\in\Lambda,\;\left( L\left( x\right) \right)_{\alpha}\circ f_{\alpha}= f_{\alpha}\circ L\left( x\right), $ then $ \left( L\left( x\right) \right)_{\alpha}\left( f_{\alpha}\left( y\right) \right) = \left( f_{\alpha}\circ L\left( x\right)\right) \left( y\right) = f_{\alpha}\left( xy\right) =f_{\alpha}\left( x\right) f_{\alpha}\left( y\right)  $ for all $ y\in E. $ Since the map $ l:\left( E_{\alpha},\overline{p}_{\alpha}\right) \rightarrow\left( M_{c}( E_{\alpha}),\lVert .\rVert_{\alpha}\right),$ \\
$l\left( x_{\alpha}\right) \left( y_{\alpha}\right) =x_{\alpha}y_{\alpha},  $ is an isometric isomorphism (into), it follows that\\
 $ \lVert\left( L\left( x\right) \right)_{\alpha}\rVert_{\alpha}= \overline{p}_{\alpha}\left( f_{\alpha}\left( x\right) \right) = p_{\alpha}\left( x\right), $ so $ q_{\alpha}\left( L\left( x\right) \right) =p_{\alpha}\left( x\right).  $
 
\noindent \textbf{} 

\noindent \textbf{Proposition 2.11.} \textit{Let $ E $ be a uniform topological algebra without unit, and let $ E_{e} $ be the algebra obtained from $ E $ by adjoining the unit. Then the map $ g:E_{e}\rightarrow M_{c}(E),\;g\left( \left( x,\lambda\right) \right) =L\left( x\right) +\lambda I $ is an algebra isomorphism (into).}

\noindent \textbf{} 

\noindent \textbf{Proof.} It is easy to show that $ g $ is an algebra homomorphism. Let $ \left( x,\lambda\right) \in E_{e} $ with $ g\left( \left( x,\lambda\right) \right) =0, $ then $ L\left( x\right) =-\lambda I. $ Suppose $ \lambda\neq 0,\;I=-\lambda^{-1}L\left( x\right) =L\left( -\lambda^{-1}x\right),$ so $ -\lambda^{-1}x $ is a left unit in $ E. $ Since $ E $ is commutative, $ -\lambda^{-1}x $ is a unit in $ E, $ a contradiction. Thus $ L\left( x\right) =0 $ and consequently $ x=0 $ since $ E $ is proper.

\noindent \textbf{} 

\noindent \textbf{Corollary 2.12.} [1, Theorem 4] \textit{Let  $\left(E,\left(p_{\alpha}\right)_{\alpha\in\Lambda}\right)$ be a uniform topological algebra without unit. Then there is a family of seminorms $ \left( s_{\alpha}\right)_{\alpha\in\Lambda} $ on $ E_{e} $ such that  $\left(E_{e},\left(s_{\alpha}\right)_{\alpha\in\Lambda}\right)$ is a uniform topological algebra and $ s_{\alpha}\left( \left( x,0\right) \right) = p_{\alpha}\left( x\right) $ for all $ x\in E $ and $ \alpha\in\Lambda. $}

\noindent \textbf{} 

\noindent \textbf{Proof.} For each $ \alpha\in\Lambda, $ we define a seminorm on $ E_{e} $ by\\
 $ s_{\alpha}\left( \left( x,\lambda\right) \right) = q_{\alpha}\left( L\left( x\right) +\lambda I\right) $ for all $ x\in E $ and $ \lambda\in\mathbb{C}. $\\
 By Propositions 2.10 and 2.11, $ \left( E_{e},\left( s_{\alpha}\right)_{\alpha\in\Lambda}\right)  $ is a uniform topological algebra and $ s_{\alpha} \left( \left( x,0\right) \right) = q_{\alpha}\left( L\left( x\right) \right) = p_{\alpha}\left( x\right) $ for all $ x\in E. $

\noindent \textbf{}

\noindent \textbf{Remark.} We have  $ s_{\alpha}\left( \left( x,\lambda\right) \right) = q_{\alpha}\left( L\left( x\right) +\lambda I\right)\leq q_{\alpha}\left( L\left( x\right) \right) +\lvert\lambda\rvert q_{\alpha}\left( I\right) = p_{\alpha}\left( x\right) +\lvert\lambda\rvert $ for all $ x\in E $ and $ \lambda\in \mathbb{C}. $ This shows that the topology on $ E_{e} $ defined by the family of seminorms $ \left( s_{\alpha}\right)_{\alpha\in\Lambda} $ is weaker than the usual topology on $ E_{e} $ defined by the family of seminorms $\left( \tilde{p}_{\alpha}\right)_{\alpha\in\Lambda}  $ where $ \tilde{p}_{\alpha}\left( \left( x,\lambda\right) \right) = p_{\alpha}\left( x\right) +\lvert\lambda\rvert. $

\noindent \textbf{}
 
\noindent \textbf{}

\noindent \textbf{References}

\noindent \textbf{} 
 
\noindent \textbf{} [1] J. Arhippainen, On locally convex square algebras, Funct. Approx. 22 (1993), 57--63.

\noindent \textbf{} [2] M. Haralampidou, The Krull nature of locally \textit{C*}-algebras, Function Spaces (Edwardsville, IL, 2002), 195--200, Contemp. Math., 328, Amer. Math. Soc., Providence, RI, 2003.

\noindent \textbf{} [3] M. Haralampidou, L. Palacios, C. Signoret, Multipliers in locally convex *-algebras, Rocky Mountain J. Math.,       43 (2013), 1931--1940.

\noindent \textbf{} [4] M. Haralampidou, L. Palacios, C. Signoret, Multipliers in  perfect locally \textit{m}-convex algebras, Banach J. Math. Anal., 9 (2015), 137--143.

\noindent \textbf{} [5] M. Haralampidou, L. Palacios, C. Signoret, Multipliers in some perfect locally \textit{m}-pseudoconvex algebras, Proceedings of the 8th International Conference on Topological Algebras and their Applications, 2014, Ed. by A. Katz, Series: De Gruyter Proceedings in Mathematics (to appear). 

\noindent \textbf{} [6] T. Husain, Multipliers of topological algebras, Dissertationes Math. (Rozprawy Mat.), 285 (1989), 36 pp.

\noindent \textbf{} [7] A. Inoue, Locally C*-algebras, Mem. Faculty Sci. Kyushu Univ. (Ser. A), 25 (1971), 197--235.

\noindent \textbf{} [8] A. Mallios, Topological Algebras. Selected Topics, North-Holland, Amsterdam, 1986.

\noindent \textbf{} [9] E. A. Michael, Locally multiplicatively convex topological algebras, Mem.
Amer. Math. Soc.,  11 (1952).

\noindent \textbf{} [10] N. C. Phillips, Inverse limits of C*-algebras, J. Operator Theory, 19 (1988), 159--195.

\noindent \textbf{} 

\noindent \textbf{} Ecole Normale Sup\'{e}rieure

\noindent \textbf{} Avenue Oued Akreuch

\noindent \textbf{} Takaddoum, BP 5118, Rabat

\noindent \textbf{} Morocco
 
\noindent \textbf{} 

\noindent \textbf{} E-mail:  mohammed.elazhari@yahoo.fr

\end{document}